\documentclass{article}
\usepackage{amssymb,amsmath}
\usepackage[mathscr]{euscript}
\input{epsf}

\newcounter{sec}

\newcounter{punct}[sec]

\def\punct{\refstepcounter{punct}{\arabic{sec}.\arabic{punct}.  }}

\newtheorem{theorem}{Theorem}[sec]
\newtheorem{proposition}[theorem]{Proposition}

\newtheorem{lemma}[theorem]{Lemma}

\newtheorem{corollary}[theorem]{Corollary}

\newtheorem{question}[theorem]{Question}
\newtheorem{conjecture}[theorem]{Conjecture}

\def\COUNTERS{\addtocounter{sec}{1}
              \setcounter{punct}{0}
          \setcounter{equation}{0}
          \setcounter{theorem}{0}
          }
          
          \def\sm{\smallskip}

\begin{document}

\newcommand{\supp}{\mathop {\mathrm {supp}}\nolimits}
\newcommand{\rk}{\mathop {\mathrm {rk}}\nolimits}
\newcommand{\Aut}{\mathop {\mathrm {Aut}}\nolimits}
\newcommand{\Out}{\mathop {\mathrm {Out}}\nolimits}
\renewcommand{\Re}{\mathop {\mathrm {Re}}\nolimits}
\newcommand{\Inn}{\mathop {\mathrm {Inn}}\nolimits}
\newcommand{\Char}{\mathop {\mathrm {Char}}\nolimits}
\newcommand{\ch}{\cosh}
\newcommand{\sh}{\sinh}
\newcommand{\Sp}{\mathop {\mathrm {Sp}}\nolimits}
\newcommand{\SOS}{\mathop {\mathrm {SO^*}}\nolimits}

\def\0{\mathbf 0}

\def\ov{\overline}
\def\wh{\widehat}
\def\wt{\widetilde}

\renewcommand{\rk}{\mathop {\mathrm {rk}}\nolimits}
\renewcommand{\Aut}{\mathop {\mathrm {Aut}}\nolimits}
\renewcommand{\Re}{\mathop {\mathrm {Re}}\nolimits}
\renewcommand{\Im}{\mathop {\mathrm {Im}}\nolimits}
\newcommand{\sgn}{\mathop {\mathrm {sgn}}\nolimits}

\def\bfa{\mathbf a}
\def\bfb{\mathbf b}
\def\bfc{\mathbf c}
\def\bfd{\mathbf d}
\def\bfe{\mathbf e}
\def\bff{\mathbf f}
\def\bfg{\mathbf g}
\def\bfh{\mathbf h}
\def\bfi{\mathbf i}
\def\bfj{\mathbf j}
\def\bfk{\mathbf k}
\def\bfl{\mathbf l}
\def\bfm{\mathbf m}
\def\bfn{\mathbf n}
\def\bfo{\mathbf o}
\def\bfp{\mathbf p}
\def\bfq{\mathbf q}
\def\bfr{\mathbf r}
\def\bfs{\mathbf s}
\def\bft{\mathbf t}
\def\bfu{\mathbf u}
\def\bfv{\mathbf v}
\def\bfw{\mathbf w}
\def\bfx{\mathbf x}
\def\bfy{\mathbf y}
\def\bfz{\mathbf z}

\def\bfA{\mathbf A}
\def\bfB{\mathbf B}
\def\bfC{\mathbf C}
\def\bfD{\mathbf D}
\def\bfE{\mathbf E}
\def\bfF{\mathbf F}
\def\bfG{\mathbf G}
\def\bfH{\mathbf H}
\def\bfI{\mathbf I}
\def\bfJ{\mathbf J}
\def\bfK{\mathbf K}
\def\bfL{\mathbf L}
\def\bfM{\mathbf M}
\def\bfN{\mathbf N}
\def\bfO{\mathbf O}
\def\bfP{\mathbf P}
\def\bfQ{\mathbf Q}
\def\bfR{\mathbf R}
\def\bfS{\mathbf S}
\def\bfT{\mathbf T}
\def\bfU{\mathbf U}
\def\bfV{\mathbf V}
\def\bfW{\mathbf W}
\def\bfX{\mathbf X}
\def\bfY{\mathbf Y}
\def\bfZ{\mathbf Z}

\def\frD{\mathfrak D}
\def\frL{\mathfrak L}
\def\frG{\mathfrak G}
\def\frg{\mathfrak g}
\def\frh{\mathfrak h}
\def\frf{\mathfrak f}
\def\frl{\mathfrak l}

\def\bfw{\mathbf w}

\def\R {{\mathbb R }}
 \def\C {{\mathbb C }}
  \def\Z{{\mathbb Z}}
  \def\H{{\mathbb H}}
\def\K{{\mathbb K}}
\def\N{{\mathbb N}}
\def\Q{{\mathbb Q}}
\def\A{{\mathbb A}}
\def\O {{\mathbb O }}

\def\T{\mathbb T}
\def\P{\mathbb P}

\def\G{\mathbb G}

\def\cD{\EuScript D}
\def\cL{\mathscr L}
\def\cK{\EuScript K}
\def\cM{\EuScript M}
\def\cN{\EuScript N}
\def\cP{\EuScript P}
\def\cQ{\EuScript Q}
\def\cR{\EuScript R}
\def\cW{\EuScript W}
\def\cY{\EuScript Y}
\def\cF{\EuScript F}
\def\cG{\EuScript G}
\def\cZ{\EuScript Z}
\def\cI{\EuScript I}
\def\cB{\EuScript B}
\def\cA{\EuScript A}

\def\bbA{\mathbb A}
\def\bbB{\mathbb B}
\def\bbD{\mathbb D}
\def\bbE{\mathbb E}
\def\bbF{\mathbb F}
\def\bbG{\mathbb G}
\def\bbI{\mathbb I}
\def\bbJ{\mathbb J}
\def\bbL{\mathbb L}
\def\bbM{\mathbb M}
\def\bbN{\mathbb N}
\def\bbO{\mathbb O}
\def\bbP{\mathbb P}
\def\bbQ{\mathbb Q}
\def\bbS{\mathbb S}
\def\bbT{\mathbb T}
\def\bbU{\mathbb U}
\def\bbV{\mathbb V}
\def\bbW{\mathbb W}
\def\bbX{\mathbb X}
\def\bbY{\mathbb Y}

\def\kappa{\varkappa}
\def\epsilon{\varepsilon}
\def\phi{\varphi}
\def\le{\leqslant}
\def\ge{\geqslant}

\def\B{\mathrm B}

\def\la{\langle}
\def\ra{\rangle}
\def\tri{\triangleright}

\def\lambdA{{\boldsymbol{\lambda}}}
\def\alphA{{\boldsymbol{\alpha}}}
\def\betA{{\boldsymbol{\beta}}}
\def\mU{{\boldsymbol{\mu}}}

\def\const{\mathrm{const}}
\def\rem{\mathrm{rem}}
\def\even{\mathrm{even}}
\def\SO{\mathrm{SO}}
\def\SL{\mathrm{SL}}
\def\PSL{\mathrm{PSL}}
\def\U{\operatorname{U}}
\def\GL{\operatorname{GL}}
\def\Mat{\operatorname{Mat}}
\def\End{\operatorname{End}}
\def\Mor{\operatorname{Mor}}
\def\Aut{\operatorname{Aut}}
\def\inv{\operatorname{inv}}
\def\red{\operatorname{red}}
\def\Ind{\operatorname{Ind}}
\def\dom{\operatorname{dom}}
\def\im{\operatorname{im}}
\def\md{\operatorname{mod\,}}
\def\indef{\operatorname{indef}}
\def\Gr{\operatorname{Gr}}

\def\arr{\rightrightarrows}

\def\cH{\EuScript{H}}
\def\cO{\EuScript{O}}
\def\cQ{\EuScript{Q}}
\def\cL{\EuScript{L}}
\def\cX{\EuScript{X}}

\def\Di{\Diamond}
\def\di{\diamond}

\def\fin{\mathrm{fin}}
\def\ThetA{\boldsymbol {\Theta}}

\def\0{\boldsymbol{0}}

\def\F{\,{\vphantom{F}}_2F_1}
\def\FF{\,{\vphantom{F}}_3F_2}
\def\H{\,\vphantom{H}^{\phantom{\star}}_2 H_2^\star}
\def\HH{\,\vphantom{H}^{\phantom{\star}}_3 H_3^\star}
\def\Ho{\,\vphantom{H}_2 H_2}

\def\disc{\mathrm{disc}}
\def\cont{\mathrm{cont}}

\def\osigma{\ov\sigma}
\def\ot{\ov t}

\begin{center}
	\bf\Large
Inner functions of matrix argument
\\ and conjugacy classes in unitary groups


\bigskip 

\sc Yury A.Neretin%
\footnote{Supported by the grant FWF  P31591.}

\end{center}

\bigskip

{\small Denote by $\B_n$ the set of complex square matrices of order $n$, whose Euclidean operator norms
are $<1$. Its Shilov boundary is the set $\U(n)$ of all unitary matrices. A holomorphic map $\B_m\to\B_n$
is inner if it sends $\U(m)$ to $\U(n)$. On the other hand we consider a group $\U(n+mj)$
and its subgroup $\U(j)$ embedded to $\U(n+mj)$ in a block-diagonal way ($m$ blocks $\U(j)$ and a unit block of size $n$). For any conjugacy class of $\U(n+mj)$
with respect to $\U(j)$ we assign a 'characteristic function', which is a rational inner map $\B_m\to\B_n$. We show that the class of inner functions, which can be obtained as
'characteristic functions', is closed with respect to natural operations
as pointwise direct sums, pointwise products, compositions, substitutions to finite-dimensional  representations
of general linear groups, etc. We also describe explicitly the   corresponding operations on conjugacy classes.}

\sm


\section{Formulation of results}

\COUNTERS

{\bf \punct Some notation.} Below:

\sm 

---  $\Mat(n)$ is the space of square complex matrices of order $n$;

\sm 

--- $1_n$ is the unit matrix of order $n$;

\sm 

--- $\GL(n,\C)$ the group of invertible complex matrices of order $n$;

\sm 

--- $\U(n)$ is the group of unitary matrices of order $n$.

\sm 

--- Let $V$, $W$ be linear spaces with bases $e_1$, \dots, $e_p$ and $f_1$, \dots, $f_q$
respectively. We order basis elements of the tensor product $V\otimes W$ as
$$
e_1\otimes f_1, \dots, e_p\otimes f_1,\,\, e_1\otimes f_2, \dots, e_p\otimes f_2,\,\, \dots,\,\,
e_1\otimes f_q, \dots, e_p\otimes f_q.
$$ 
According this we write tensor products of matrices.

\sm 

{\bf \punct Matrix balls and inner functions.}
Denote by $\|\cdot\|$ the operator norm in Euclidean space, i.e., $\|z\|^2$
is the maximal eigenvalue of the matrix $z^*z$.
Denote $\B_n$ (a {\it matrix ball}) the set of complex matrices $z$ of order $n$ such that $\|z\|<1$, by $\ov\B_n$
its closure, i.e., the set of matrices satisfying $\|z\|\le 1$.
By $\partial \B_n$ we denote the boundary of $\B_n$, i.e. the set of matrices with norm 1. 
The unitary group $\U(n)$ is contained in $\partial\B_n$ and is the Shilov boundary
of $\B_n$.

Recall that  the pseudounitary group  $\U(n,n)$   acts on $\B_n$ by biholomorphic 
transformations, and the space $\B_n$ is the  symmetric space
$$
\B_n\simeq \U(n,n)\bigr/(\U(n)\times \U(n)),
$$
(see, e.g., \cite{P-Sh}, Sect.6, \cite{Ner-book}, Sect. 2.3, see also below Subsect. \ref{ss:lin-fract}).

We say that a holomorphic map $F:\B_m\to\B_\alpha$ is {\it inner} if
its limit values on $\U(m)$ are defined a.s. and $F$ sends $\U(m)$ to $\U(\alpha)$.
Below we discuss only rational maps, so a meaning of the term 'limit values' here is clear. 

\sm 

{\sc Remark.}
Recall that inner functions $\B_1\to \B_1$ (i.e., holomorphic maps of the
unit disk $|z|<1$ to itself that also send the circle $|z|=1$ to itself) are a classical topic of function theory of complex 
variable, see, e.g. \cite{Gar}. Inner functions $\B_1\to \B_\alpha$ arose in the context of works
of M.S.Livshits 1946-1954 on spectral theory of operators closed to unitary operators, see \cite{Liv1},
\cite{Liv2}, see also \cite{Nad-F}. V.P.Potapov \cite{Pot} obtained a multiplicative
representation of such functions, see also \cite{Alp}. Inner functions $\B_m\to \B_\alpha$
arose in \cite{Ner-AAMP}--\cite{Ner-spher} in representation theory of infinite-dimensional
classical groups.
\hfill $\boxtimes$

\sm 

{\bf\punct Colligations and characterictic functions.}
 Fix $\alpha$ and $m\in \N$.  $j=0$, $1$, $2$, \dots.
Consider a unitary group $\U(\alpha+mj)$ and its subgroup $\U(j)$ embedded
as
$$
T\mapsto \begin{pmatrix}
1_\alpha&0\\0& T\otimes 1_m
\end{pmatrix}:=
\begin{pmatrix}
1_\alpha&0&0&\dots&0\\
0&T&0&\dots&0\\
0&0&T&\dots&0\\
\vdots&\vdots&\vdots&\ddots&\vdots\\
0&0&0&\dots&T
\end{pmatrix}\in \U(\alpha+mj).
$$
Consider conjugacy classes of the group $\U(\alpha+mj)$ with respect to the subgroup $\U(j)$, i.e., matrices
defined up to the equivalence
$$
g\sim hgh^{-1}, \qquad \text{where $g\in \U(\alpha+mj)$, $h\in \U(j)$.}
$$
We call such conjugacy classes by {\it colligations}.
Let $S\in \B_m$, let $s_{\mu\nu}$ be its matrix elements. We write $g$
as a block matrix of size $\alpha+\underbrace{j+\dots+j}_{\text{$m$ times}}$,
\begin{equation}
g=\begin{pmatrix}
a&b_1&\dots\\c_1&d_{11}&\dots\\\vdots&\vdots&\ddots
\end{pmatrix}\in \U(\alpha+mj),
\label{eq:gabsd}
\end{equation}
and consider the following relation:
\begin{equation}
\begin{pmatrix}
p\\ x_1\\\vdots\\ x_m
\end{pmatrix}
=
\begin{pmatrix}
a&b_1&\dots&b_m\\
c_1&d_{11}&\dots&d_{1m}\\
\vdots&\vdots&\ddots&\vdots\\
c_k&d_{m1}&\dots&d_{mm}
\end{pmatrix}
\begin{pmatrix}
q\\
s_{11}x_1+\dots +s_{1m}x_m
\\
\vdots \\
s_{m1}x_1+\dots +s_{mm}x_m
\end{pmatrix},
\label{eq:basic}
\end{equation}
where columns $p$, $q\in \C^\alpha$, and $x_1$, \dots, $x_m\in \C^j$.
We eliminate variables $x_1$, \dots, $x_m$ and get a dependence
$$
p =\Theta[g;S] q,
$$
where $\Theta[g;S]$ is a rational matrix-valued function of the variable $S$ depending on a parameter $g$.
By \cite{Ner-AAMP}, Theorem 4.1, this function depends only on the conjugacy class containing $g$  and
is an inner function of the matrix variable $S$. We call $\Theta[g;S]$ by a {\it characteristic function of a colligation.}

Let us repeat the definition in other terms. Denote
$$
1_j\otimes S:=
\begin{pmatrix}
s_{11}\cdot 1_j&\dots&s_{1m}\cdot 1_j\\
\vdots&\ddots&\vdots\\
s_{m1}\cdot 1_j&\dots&s_{mm}\cdot 1_j
\end{pmatrix}.
$$
Then
\begin{equation}
\Theta[g;S]=a+b \,(1_j\otimes S)\,\Bigl(1_{mj}-d\, (1_j\otimes S) \Bigr)^{-1} c.
\label{eq:chi-otimes}
\end{equation}



\sm

{\sc Remark.} Let $H=H_1\oplus H_2$ be a Hilbert space,
let $\begin{pmatrix}
a&b\\c&d
\end{pmatrix}$ be a  unitary operator in $H$. 
 The {\it Livshits characteristic function},  which arose in spectral theory of
 	non-normal 
 operators, see \cite{Liv1}, \cite{Liv2} (see also \cite{Nad-F}, \cite{Bro1}, \cite{Bro2})
 is given by 
 \begin{equation}
\Theta(\lambda)=a+b(1-\lambda d)^{-1}c.
\label{eq:char}
\end{equation}
In our construction this corresponds to $m=1$. We prefer to use
the original term 'characteristic function', which  emphasize
analogy with characteristic numbers and characteristic polynomials.
Our equation \eqref{eq:basic} can be regarded as an extension
of the equation $Ax=s x$. But the term 'characteristic function' is overburdened (it has
two other common meanings: indicator functions and Fourier transforms of measures in probability).
  There is another term {\it transfer function} (see, e.g., \cite{Bart})
for \eqref{eq:char}, which came from system theory.
\hfill $\boxtimes$

\sm

{\sc Remark.} Colligation type structures arise in representation theory of infinite dimensional
classical groups as 'semigroups of double cosets'. Firstly, this was observed by G.~I.~Olshanski
in \cite{Olsh}, see also \cite{Ner-cat}, Sect.IX.3-4. Such semigroups act in spaces of unitary representations of corresponding classical
groups by certain operators with  Gaussian kernels. A Gaussian kernel is determined by a matrix, and Olshanski
showed that such matrices are given by expressions similar to matrix-valued characteristic functions
of one variable. In fact, the origin of inner functions of matrix variables
in \cite{Ner-AAMP}, \cite{Ner-spher} was similar, but an initial point was a more general class of unitary representations
from  N.~I.~Nessonov \cite{Ness}. From this point of view basic objects
are semigroups of 'colligations' and inner functions are a tool for  understanding of colligations.
 In considerations of this paper we do not refer to representation theory and regard 
inner function of matrix variables and colligations  as  abstract topics.
\hfill $\boxtimes$

\sm

{\bf \punct A conjecture.} Denote by 
$$\Inn(m,\alpha)=\Inn[\B_m,\B_\alpha]$$
 the space of all rational
interior maps $F:\ov\B_m\to\ov\B_\alpha$. Denote by $\Inn_\circ(m,\alpha)$
its subset consisting of maps $F$ such that $F(\B_m)\subset \B_\alpha$.
By the maximum modulus principle the last condition is equivalent to: for some $z_0\in\B_m$
we have $F(z_0)\in \B_\alpha$.

We also define the space 
$$\Char(m,\alpha)=\Char[\B_m,\B_\alpha]$$
 consisting of characteristic functions
determined by all possible elements of $\U(\alpha+mj)$ for $j=0$, 1, 2, \dots.
By 
$\Char_\circ(m,\alpha)$
 we denote its subset, consisting of function
$\Theta[g;\cdot]$ sending $\B_m\to\B_\alpha$. In notation
\eqref{eq:gabsd}, a map $\Theta[g;\cdot]$ is contained $\Char_\circ(m,\alpha)$
if and only if $\|a\|<1$.

\begin{conjecture}
	Any rational inner function is a characteristic function of some colligation, i.e.,
$$\Inn(m,\alpha)=\Char(m,\alpha).$$	
\end{conjecture}	

{\sc Remark.}
The conjecture was formulated in \cite{Ner-spher}. It is not doubtless
since similar statement for inner functions in polydisks is false
(or is valid under additional conditions to rational inner functions, cf. the case 
of polysisk in \cite{Kne}, \cite{BB}).
\hfill $\boxtimes$

\sm

{\bf \punct Operations in $\Char(m,\alpha)$.}
In this paper we chow that the class $\coprod_{m,\alpha}\Char(m,\alpha)\subset \coprod_{m,\alpha}\Inn(m,\alpha)$
is closed with respect to several natural operations%
\footnote{These statements were announced in \cite{Ner-spher}.}.
In all cases {\it we describe explicitly operations over colligations corresponding to operations
over inner functions,} these formulas are presented in proofs. 

\begin{theorem}
	\label{th:sums}
	{\rm a)} Let $F_1\in \Char(m,\alpha)$, $F_2\in \Char(m,\beta)$. Then
	$F_1\oplus F_2\in \Char(m,\alpha+\beta)$.
	
	\sm
	
	{\rm b)} Let $F\in \Char(m,\alpha+\beta)$ admits a decomposition into a direct sum
	$F=F_1\oplus F_2$, where $F_1\in \Inn(m,\alpha)$, $F_2\in \Inn(m,\beta)$.
	Then $F_1\in \Char(m,\alpha)$, $F_2\in \Char(m,\beta)$.
\end{theorem} 

The first (trivial) part of the statement is proved in Subsect. \ref{ss:sums-1}, the second part
in Subsect. \ref{ss:sums-2}.
 
The following statement was obtained in \cite{Ner-AAMP}, Theorem 4.1.

\begin{theorem}
	\label{th:products}
	Let $F_1$, $F_2\in \Char(m,\alpha)$. Then the pointwise  product $F_1 F_2$ of matrix-valued 
	functions $F_1$, $F_2$
	 is contained
	in $\Char(m,\alpha)$. 
\end{theorem}

\begin{theorem}
	\label{th:tensors}
Let $F_1\in \Char(m,\alpha)$, $F_2\in \Char(m,\beta)$. Then the pointwise tensor product $F_1\otimes F_2$ of matrix-valued 
functions $F_1$, $F_2$
is contained
in $\Char(m,\alpha\beta)$. 
\end{theorem}

This is proved in Subsect. \ref{ss:tensors}.

\begin{theorem}
	\label{th:compositions}
	Let $G\in \Char( \beta,\gamma)$ be defined by a matrix $\begin{pmatrix}
	a&b\\c&d
	\end{pmatrix}\in \U(\gamma+\beta j)$,  $F\in \Char(\alpha,\beta)$ by a matrix
	$\begin{pmatrix}
	p&q\\r&t
	\end{pmatrix}\in \U(\beta+\alpha i)$. Let 
	\begin{equation}
	\label{eq:cond1}
	\det\bigl(1_{\beta j}-d\,(1_j\otimes p)\bigr)\ne 0.
	\end{equation} 
	Then $G\circ F\in \Char(\alpha,\gamma)$. More generally, the conclusion holds
	if 
	\begin{equation}
	\label{eq:cond2}
	\det\bigl(1_{\beta j}-d\, (1_j\otimes F(S_0)) \bigr)\ne 0\qquad\text{for some $S_0\in \B_\alpha$.}
	\end{equation}
\end{theorem}

{\sc Remarks.} a) In particular, the condition \eqref{eq:cond1}  holds if $\|d\|<1$ or if $\|p\|<1$.
Recall that $\|d\|\le 1$, $\|p\|\le 1$.

\sm

b) It can happened that the image of $F$ is contained in the set of discontinuity
of $G$. The condition \eqref{eq:cond2} is sufficient (and not necessary) for avoiding this situation. 
\hfill $\boxtimes$

\sm

Theorem is proved in Subsect. \ref{ss:compositions}.

\sm 

Next, consider a unitary finite-dimensional representation $\rho$ of a unitary group $\U(n)$.  Then (see, e.g., \cite{Zhe}, \S42) it admits a unique holomorphic continuation to a representation of $\GL(n,\C)$. The representation $\rho$ is called {\it polynomial}, if all  matrix 
elements of $\rho(g)$ are polynomials in matrix elements of $g\in\GL(n,\C)$. Consider the semigroup
$\Mat^\times(n)$ of all matrices of order $n$ with respect to the multiplication. The group
$\GL(n,\C)$ is dense in $\Mat^\times(n)$ and all polynomial representations of $\GL(n,\C)$
have a continuous extensions to the semigroup $\Mat^\times(n)$ (matrices $\rho(\cdot)$ are determined
by the same polynomials).

\begin{theorem}
	\label{th:representations}
	Let $\rho$ be a polynomial unitary  representation of $\U(\alpha)$. Let $F\in\Inn(m,\alpha)$.
	Then $\rho\circ F$ is contained in $\Inn(m,\dim\rho)$.
\end{theorem}

The statement is proved in Subsect. \ref{ss:tensors-over}.

\begin{corollary}
	Let $F\in\Inn(m,\alpha)$. Then $\det(F)\in  \Inn(m,1)$.
\end{corollary}

{\bf \punct Some remarks on behavior of inner functions on strata of boundaries.}
Recall (see \cite{P-Sh}, Sect.6, see below Subsect. \ref{ss:boundary}) that the boundary of the domain $\B_m\subset \Mat(m)$ is a disjoint union of a continual family
of  complex (open) manifolds (boundary components), these components are maximal
complex manifolds that are contained in the boundary. Each component $C$ is 
biholomorphically equivalent to some matrix ball $\B_\nu$, where $\nu=0$, $1$, \dots $m-1$,
and a biholomorphic map $\B_\mu\to C$ extends continuously to a homeomorphism of closures 
$\ov \B_\mu\to \ov C$.

The following statements are obvious.

\begin{proposition}
	{\rm a)} Let $F\in \Inn[\B_m,\B_\alpha]\setminus \Inn_\circ[\B_m,\B_\alpha]$. Then
	$F(\B_m)$ is contained in a unique boundary component $C\subset \B_\alpha$.
	Moreover, $F$ is contained in $\Inn_\circ[B_m,C]$.
	
	\sm 
	
	{\rm b)} Let $F\in \Inn[\B_m,\B_\alpha]$ be continuous at some point of a boundary component $C\subset \B_m$.
	Then $F\in \Inn[C,\B_\alpha]$.
\end{proposition}

The next statement is proved in Subsect. \ref{ss:boundary-restriction}.

\begin{theorem}
	\label{th:boundary}
	{\rm a)} Let $F\in \Char[\B_m,\B_\alpha]$ send $\B_m$ into a boundary component $C\subset \ov\B_\alpha$.
	Then $F\in \Char_\circ[\B_m,C]$.
	
	\sm 
	
	{\rm b)} Let $F\in \Char[\B_m,\B_\alpha]$ be determined by a matrix $\begin{pmatrix}
	a&b\\c&d
	\end{pmatrix}\in \U(\alpha+mj)$. Let $C\subset \ov\B_m$ be a boundary component.
	 Let $\det\bigl(1_{mj}- d(1_j\otimes S_0)\bigr)\ne 0$ for some $S_0\in C$. Then
	 the restriction of $F$ to $C$ is contained in $\Char[C,\B_\alpha]$.
\end{theorem}

\sm

{\bf\punct On some extensions of the construction.}
As we mentioned above, $\B_n$ is a symmetric space. Our construction of inner functions can be automatically
extended to Hermitian symmetric spaces of series
$$
\U(p,q)/\bigr(\U(p)\times \U(q)\bigl),
\qquad 
\Sp(2n,\R)/\U(n),\qquad \SOS(2n,\R)
$$ 
(see, e.g., \cite{P-Sh}, they are called classical complex domains of I, II, and III types)
and for direct products of such spaces (see  \cite{Ner-AAMP}, \cite{Ner-char}).

However, initial spaces and target spaces in such constructions are not independent.
For instance,  this approach does not produce inner functions from 
the usual unit ball $\U(n,1)/\bigl(\U(n)\times \U(1)\bigr)$ to the unit disk (they exist
according A.B.Alexandrov \cite{Ale} and E.L\o w \cite{Low}), and, more generally,
for functions from $\U(n,1)/\bigl(\U(n)\times \U(1)\bigr)$ to the unit disk for $p\ne q$ (they also exist
according Alexandrov \cite{Ale2}). 

However it is possible to take $m=\infty$.
There arises the following question.

\begin{question}
	Let $g$ be a unitary operator  
	$$\C\oplus \ell^2\oplus \ell^2\,\to \C\oplus \ell^2,$$
	i.e., $gg^*=1$, $g^*g=1$. Represent $g$ in the block form
	$$
	g=\left(
\begin{array}{c|cc}
a&b_1&b_2\\
\hline
c&d_1&d_2
\end{array}\right)
	.$$
Let $(s_1,s_2)$ be a point of the open unit ball in $\C^2$. Set
\begin{multline}
\Theta\bigl[g;(s_1,s_2)\bigr]:=\\:=
a+\begin{pmatrix}  b_1& b_2
\end{pmatrix}
\begin{pmatrix}
s_1\cdot 1_\infty\\s_2\cdot 1_\infty
\end{pmatrix}
+
\left(1_\infty- \begin{pmatrix}  d_1& d_2
\end{pmatrix}
\begin{pmatrix}
s_1\cdot 1_\infty\\s_2\cdot 1_\infty
\end{pmatrix} \right)^{-1}c=\\
a+\begin{pmatrix} s_1 b_1+s_2 b_2
\end{pmatrix}
\bigl[1- s_1d_1-s_2d_2\bigr]^{-1}c.
\label{eq:ball-inner}
\end{multline}	
Is it possible to find inner function of such a type?
Is it possible to find conditions for $g$ under which $\Theta[g]$ is an inner function
in the unit ball in $\C^2$?
\end{question}

An argument for this supposal is very simple. We write the following relation
$$
\begin{pmatrix}
p\\x
\end{pmatrix} =\left(
\begin{array}{c|cc}
a&b_1&b_2\\
\hline
c&d_1&d_2
\end{array}\right)
\begin{pmatrix}
q\\s_1 x\\s_2 x
\end{pmatrix}.
$$
Eliminating variables $x$, we come to $p=\Theta\bigl[g,(s_1,s_2)\bigr]\,q$.
On the other hand, the matrix $g$ is unitary, therefore
$|p|^2+\|x\|^2=|q|^2+|s_1|^2 \|x\|^2+ |s_2|^2 \|x\|^2$, i.e.,
$$ 
|p|^2=|q|^2-(1-|s_1|^2-|s_2|^2)\|x\|^2.
$$
If $|s_1|^2+|s_2|^2<1$, then $|p|\le |q|$ and $\bigl|\Theta\bigl[g,(s_1,s_2)\bigr]\bigr|\le1$.
At first glance it seems that $|s_1|^2+|s_2|^2=1$ immediately implies $|p|= |q|$.
But it is not that simple, since the matrix in square brackets in \eqref{eq:ball-inner}
can be non-invertible in this case.

It is easy to present examples, when  function   $\Theta[g;\cdot]$
is not inner, however a Livshits characteristic function
 \eqref{eq:char} also not always is inner, see \cite{Nad-F}, \S VI.1.

\section{Preliminaries}

\COUNTERS

{\bf \punct Linear fractional maps.%
\label{ss:lin-fract}} We realize the {\it pseudo-unitary group}
$\U(n,n)$ as the group of all complex block matrices
$g=\begin{pmatrix}
A&B\\C&D
\end{pmatrix}$ of size $(n+n)$ satisfying the condition
$$
g \begin{pmatrix}
-1_n&0\\0&1_n
\end{pmatrix}g^*=\begin{pmatrix}
-1_n&0\\0&1_n
\end{pmatrix}.
$$
For each $g\in\U(n,n)$ we consider the following linear fractional transformation
of the space $\Mat(n)$:
\begin{equation}
\gamma[g;z]=(A+zC)^{-1}(B+zD).
\label{eq:lin-fra}
\end{equation}
Such transformations send the matrix ball $\B_n$ to itself
(see, e.g., \cite{P-Sh}, \S6, \cite{Ner-book}, Sect. 2.3). This action of $\U(n,n)$ is transitive, the stabilizer of the point $z=0$
consists of matrices $\begin{pmatrix}
a&0\\0&d
\end{pmatrix}$, where $a\in\U(n)$, $d\in \U(n)$.
So $\B_n$ is a homogeneous space 
$$
\B_n=\U(n,n)/\bigl(\U(n)\times\U(n)\bigr).
$$

{\bf \punct Realization of $\B_n$ as a domain in a Grassmannian.}
Consider the pseudo-Euclidean space 
$$V^{2n}=V_-^n\oplus V_+^n:=\C^n\oplus \C^n$$
 equipped with the Hermitian form $\cM=\cM_n$ determined by the matrix
 $\begin{pmatrix}
 -1_n&0\\0&1_n
 \end{pmatrix}$. Our group $\U(n,n)$ preserves this form.

 Denote $\Gr(n)$ the Grassmanian  of all $n$-dimensional subspaces in $V^{2n}$.
 We say that a subspace $L\in\Gr(n)$ is {\it negative} ({\it semi-negative}) if the form $\cM$ is 
 negative (resp., semi-negative) defined on $L$. A subspace $L\in\Gr(n)$ is {\it isotropic} 
 if  the form $\cM$ is zero on $L$. We denote the corresponding subsets in
 the Grassmannian by
 $$\Gr^{<0}(n), \qquad \Gr^{\le 0}(n), \qquad \Gr^{0}(n)$$
 respectively.
 
  For any linear map
 $z:V_-^n\to V_+^n$ we consider the $n$-dimensional space $L[z]$ consisting of vectors of the
 form $(v_-,v_-z)$. If $z\in \B_n$, then the form $\cM$ is  negative
 on $L(z)$. Vice versa, any $n$-dimensional negative subspace in $V$ has the form
 $L(z)$ for some $z\in \B_n$. Formula \eqref{eq:lin-fra} corresponds to the natural
 action of $\U(n,n)$ on the set of negative subspaces. Also
 $$
 \ov \B_n\simeq \Gr^{\le 0}(n),\qquad \U(n)\simeq \Gr^0(n).
 $$
 
 \sm 

{\bf \punct The structure of the boundary of the matrix ball $\ov\B_n$.%
\label{ss:boundary}} The boundary of $\B_n$ consists
of $n$ orbits $\cO_j$, where $j=1$, 2, \dots, $n-1$, of the group $\U(n,n)$, representatives of orbits are matrices
of the form $\begin{pmatrix}
1_j&0\\
0& 0_{n-j}
\end{pmatrix}$. The Shilov boundary $\U(n)$ corresponds to $j=n$.

On the language of the Grassmannian $\Gr(n)$ 
the orbit $\cO_j$ corresponds to semi-negative subspaces $L$ such that
the rank the form $\cM$ on $L$ is $(n-j)$.

Any component of the boundary $\partial \B_n$ can be reduced by a linear fractional transformation
 to the form
$$
\begin{pmatrix}
u&0\\0& 1_j
\end{pmatrix}, \qquad\text{where $u$ ranges in $\B_{n-j}$}.
$$

On the language of the Grassmannian, boundary components $C$ are enumerated
by a number $j$ and
 a $j$-dimensional isotropic subspace $W\subset \C^n\oplus \C^n$.
 The corresponding component consists of
 all $n$-dimensional $\cM$-semi-negative subspaces $L\supset W$ such that
 the kernel of $\cM$ on $L$ coincides with $W$.

\sm 

{\bf\punct Linear relations.}
Let $V$, $W$ be linear spaces. A {\it linear relation}
$Y:V\arr W$ is a linear subspace in $V\oplus W$. For 
linear relations $Y:V\arr W$ and $Z:W\arr U$
we define their {\it product} $ZY:V\arr U$
as the set of all $v\oplus u\in V\oplus U$, for which
there exists $w\in W$ satisfying $v\oplus w\in Y$, $w\oplus u\in Q$.

Let $H\subset V$ be a linear subspace. The subspace $YH\subset W$ consists
of $w\in W$, for which there exists $v\in H$ satisfying $v\oplus w\in Y$. We can consider
$H$ as linear relation $0\arr V$, therefore we can understand $PH$ as a product of linear relations.

For a linear relation $Y:V\arr W$ we define

\sm 

--- {\it kernel} $\ker Y\subset V$ as the intersection $Y\cap V$;

\sm 

--- {\it domain} $\dom Y\subset V$ is the image of the projection $Y$ to $V$ along $W$;

\sm 

--- {\it image} $\im Y\subset W$ is the image of projection of $Y$ to $W$;

\sm 

--- {\it indefiniteness} $\indef Y$ is $Y\cap W$.

\sm

A product of linear relations $Y:V\arr W$, $Z:W\arr Y$
is a continuous operation $(Y,Z)\to ZY$ outside sets 
$$\ker Z\cap \indef Y\ne 0\qquad \text{and }\qquad \im Y+\dom Z\ne W.$$

\sm 

{\bf\punct Isotropic category.}
Objects of the {\it isotropic category}, see \cite{Ner-book}, Sect. 2.10, are spaces 
$$V^{2n}=V_+^n\oplus V_-^n\simeq \C^n\oplus\C^n,$$
where $n=0$, 1, 2, \dots. 
A morphism $V^{2n}\to V^{2m}$ is a linear relation $Y:V^{2n}\arr V^{2m}$
satisfying conditions:

\sm 

1) if $v\oplus v'\in Y$, then $\cM_n(v,v)=\cM_m(v',v')$;

\sm 

2) $\dim P$ is maximal possible, i.e., $\dim P=m+n$.

\sm 

A product of morphism is the product of linear relations.
The group of automorphisms of $V^{2n}$ is $\U(n,n)$.

Emphasize that the product has points of discontinuity.

Equip $V^{2n}\oplus V^{2m}$ with the difference of Hermitian forms in this space,
$$
\cM_{m,n}(v\oplus w, v'\oplus w')
:=\cM_n(v,v')-\cM_m(w,w').
$$
Then the subspace $V_-^n\oplus V_+^m\subset V^{2n}\oplus V^{2m}$ is  negative with respect
to the form $\cM_{m,n}$ and the subspace $V_+^n\oplus V_-^m\subset V^{2n}\oplus V^{2m}$
is positive. So we can apply above reasoning and get:

\sm

{\it A relation $P:V^{2n}\arr V^{2m}$ is isotropic if and only if $P$ is a graph
	of a unitary operator  $V_-^n\oplus V_+^m\to V_+^n\oplus V_-^m$.}

\sm 

Thus the set of morphisms from $V^{2n}$ to $V^{2m}$ is in one-to-one correspondence 
with the unitary group $\U(n+m)$ and the product of morphisms $Y:V^{2n}\arr V^{2m}$,
$Z:V^{2m}\arr V^{2k}$ induces an operation
$$
\U(n+m)\times \U(m+k)\to \U(n+k).
$$

For the following statement, see, e.g., \cite{Ner-book}, Theorem 2.8.4.

\begin{proposition}
	\label{pr:isotropic-product}
 Let 
$Y:V^{2k}\arr V^{2m}$ corresponds to a unitary matrix $\upsilon=\begin{pmatrix}
p&q\\r&t
\end{pmatrix}\in \U(k+m)$ and $Z:V^{2m}\arr V^{2n}$ correspond to
a unitary matrix $\zeta=\begin{pmatrix}a&b\\c&d\end{pmatrix}\in \U(n+m)$.
Let
$$
\det(1-pd)^{-1}\ne 0.
$$
Then $Z Y $ corresponds to the matrix
\begin{equation}
\zeta\circledast\upsilon
:=\begin{pmatrix}
a+b(1-pd)^{-1}pc& b(1-pd)^{-1} q\\
r(1-dp)^{-1}c&t+rd(1-pd)^{-1}q
\end{pmatrix}.
\label{eq:circledast}
\end{equation}
\end{proposition}

\sm 

{\bf \punct Krein--Shmul'yan maps.}
Let $L\in\Gr^{\le 0}(m)$. Applying a morphism
$Z:V^{2m}\arr V^{2n}$ of the isotropic category to $L$ we again get an element
of $\Gr^{\le 0}(n)$, so we get a map $\ov\B_m\to \ov\B_m$ (see \cite{Ner-book}, Theorem 2.9.1).
Let  $\zeta=\begin{pmatrix}
a&b\\c&d
\end{pmatrix}
$ be the unitary matrix corresponding to $Z$. Then {\it the corresponding
map $\sigma[\zeta]$ is given by the formula
\begin{equation}
\sigma[\zeta;u]:\, u\mapsto a+bu(1_m-ud)^{-1}c, \qquad\text{where $u\in \ov\B_m$}.
\label{eq:kr-sh}
\end{equation}
This holds if $u$ satisfies the condition $\det(1-ud)\ne 0$.}

Notice, that

\sm 

1) for any $\zeta\in\U(n+m)$
our map is continuous as a map $\B_m\to\ov\B_n$;

\sm 

2) if $\|d\|<1$, the our formula determines a continuous map 
$\ov\B_m\to\ov\B_n$;

\sm 

3) if $\|a\|<1$, then our formula determines a continuous map
$\B_m\to\B_n$.

\sm

{\sc Remark.} A map $\sigma[\zeta;z]$ is a special case of
Krein--Shmul'yan maps, see \cite{KSh}, see also \cite{Ner-book}, Sect. 2.9.
\hfill $\boxtimes$

\sm

{\sc Remark.} A map \eqref{eq:kr-sh} $\B_m\to \B_m$ is inner and moreover
its characteristic function is determined by an element of $\U(n+m\cdot 1)$.
\hfill $\boxtimes$

\sm


\begin{lemma}
	\label{l:krein-schm}
	Let $\zeta$ and $\upsilon$ be the same as in Proposition \ref{pr:isotropic-product}.
	Then for any $u\in \B_k$ we have
	\begin{equation}
	\sigma\bigl[\zeta;\,\sigma[\upsilon;\,u]  \bigr]=
	\sigma\bigl[\zeta\circledast \upsilon;\, u  \bigr].
	\label{eq:sigma-sigma}
	\end{equation}
\end{lemma}

{\sc Remark.} Cf. \cite{Ner-book}, Theorem 2.9.4, but conditions of this theorem are not satisfied.
Basically, Lemma \ref{l:krein-schm} claims  the associativity of product of linear relations
$0\arr V^{2k}\arr V^{2m}\arr V^{2n}$. However, formula \eqref{eq:circledast} is not valid 
on the surface $\det(1-pd)=0$ and this requires some care. To avoid references to proofs
or repetitions of proofs, we present a formal calculation. 
\hfill $\boxtimes$

\sm

{\sc Proof.} We must transform the following expression to the Krein--Shmul'yan form:
\begin{equation}
a+b\, \Bigl(u(1-tu)^{-1} \Bigr|_{u=p+qz(1-tz)^{-1}r} \Bigr)\, c.
\label{eq:subs}
\end{equation}

{\sc Step 1.}
It is sufficient to examine the expression in big brackets, it is a sum $I+J$ of two summands
\begin{align*}
I:&=p \Bigl[1-dq-dq z(1-tz)^{-1}r\Bigr]^{-1};
\\
J:&=qz(1-tz)^{-1} \Bigl[1-dq-dq z(1-tz)^{-1}r\Bigr]^{-1}.
\end{align*}
First, we must show that the inverse matrix $\bigl[\dots\bigr]^{-1}$ exists.
Since $(1-dp)^{-1}$ is invertible, we can transform $\bigl[\dots\bigr]^{-1}$
as
\begin{equation}
\bigl[\dots\bigr]^{-1}=(1-dq)^{-1}\Bigl(1- dq z(1-tz)^{-1}\cdot r(1-dq)^{-1}\Bigr)^{-1}
\label{eq:skobka}
\end{equation}
Next, we notice that matrices $(1-AB)$ and $(1-BA)$ are invertible or noninvertible simultaneously.
Therefore it is sufficient to verify existence of the matrix
$$
\Bigl(1-  r(1-dq)^{-1}\cdot dq z(1-tz)^{-1}\Bigr)^{-1}
=(1-tz) \Bigl(1- \bigl\{t+r(1-dq)^{-1}\cdot dq\bigr\}\cdot z\Bigr)^{-1}.
$$
Since $\|t\|\le 1$, $\|z\|<1$ the matrix $(1-tz)$ is invertible.
Next, the expression in curly brackets coincides with lower right block of the matrix
\eqref{eq:circledast}. Therefore $\|\{\dots \}\|\le 1$ and the second factor is well-defined.

\sm

{\sc Step 2.}
Transforming the summand $J$ with \eqref{eq:skobka} we get
$$
J=q z(1-tz)^{-1}\cdot r(1-dq)^{-1}\Bigl(1- dq z(1-tz)^{-1}\cdot r(1-dq)^{-1}\Bigr)^{-1}.
$$
Keeping in mind the matrix identity
\begin{equation*}
A(1-BA)^{-1}=(1-AB)^{-1}A,
\end{equation*}
we come to
\begin{multline*}
J=q z(1-tz)^{-1}\cdot \Bigl(1- r(1-dq)^{-1}\cdot dq z(1-tz)^{-1}\Bigr)^{-1}r(1-dq)^{-1}
=\\=qz \Bigl(1- \bigl\{t+r(1-dq)^{-1}\cdot dq\bigr\}\cdot z\Bigr)^{-1}r(1-dq)^{-1}.
\end{multline*}
Next, we transform the summand $I$ with \eqref{eq:skobka}
and apply the identity
\begin{equation}
(1-C)^{-1}=1+C(1-C)^{-1}
\label{eq:CCC}
\end{equation}
to the second factor in \eqref{eq:skobka}.
We come to
\begin{multline*}
I=
p(1-dp)^{-1}+\\+
 p(1-dp)^{-1} d\cdot q z(1-tz)^{-1} r(1-dq)^{-1} \Bigl(1- dq z(1-tz)^{-1} r(1-dq)^{-1}\Bigr)^{-1}
=\\=p(1-dp)^{-1}+ p(1-dp)^{-1} d\cdot J
\end{multline*}
Thus
\begin{multline*}
I+J=p(1-dp)^{-1}+\Bigl\{p(1-pd)^{-1}d+1\Bigr\}\cdot J=\\=
(1-pd)^{-1}p+\bigl\{(1-pd)^{-1}\bigr\}\cdot qz \Bigl(1- \bigl\{t+r(1-dq)^{-1}\cdot dq\bigr\}\cdot z\Bigr)^{-1}r(1-dq)^{-1}
\end{multline*}
(we applied \eqref{eq:CCC} to the expression in curly brackets).
We substitute the result to  \eqref{eq:subs} instead of the expression in big brackets
and get the desired formula.
\hfill $\square$

\sm








{\bf \punct Characteristic functions and Krein--Shmul'yan maps.}
The formula \eqref{eq:chi-otimes} can be written as
\begin{equation}
\Theta\left[\begin{pmatrix}
a&b\\c&d
\end{pmatrix};S \right]=
\sigma\left[\begin{pmatrix}
a&b\\c&d
\end{pmatrix};1_j\otimes S \right].
\label{eq:chi-sigma}
\end{equation}

{\bf \punct Polynomial representations of $\GL(n,\C)$.%
\label{ss:representation}}
Irreducible holomorphic representations $\rho_\bfm(g)$ of $\GL(n,\C)$ are enumerated
by 'signatures' 
$$
\bfm:=(m_1, \dots, m_n), \quad \text{where $m_j\in \Z$ and $m_1\ge m_2\ge\dots\ge m_n$,}
$$
see, e.g., \cite{Zhe}, \S49-50. 
Recall a semi-explicit construction of $\rho_\bfm$.

Consider the space $\C^n$
with the standard basis $e_1$, \dots, $e_n$.  
The representation $\lambda_k(g)$  with signature
$(\underbrace{1,\dots,1}_{\text{$k$ times}}, 0,\dots,0)$ 
is the representation in the $k$-th exterior power $\bigwedge^k \C^n$, its highest weight
vector is $v_k=e_1\wedge \dots \wedge e_k$. The matrix element
$\la \rho(g) v_k,v_k\ra$ is the $k$-th principal minor $\Delta_k(g)$ of $g$.
The representation $\lambda_n$ is simply $\det(g)$ (in particular, we can consider its negative tensor powers). 

The representation $\rho_\bfm$ is a subrepresentation of
\begin{equation}
\bigotimes_{k=1}^{n-1} \lambda_k^{\otimes (m_k-m_{k+1})}(g)\,\otimes \,\lambda_k(g)^{m_n}. 
\label{eq:bigotimes}
\end{equation}
More precisely, $\rho_\bfm$ is the cyclic span of the highest weight vector
$$
\xi_\bfm:=\bigotimes_{k=1}^{n-1} (e_1\wedge \dots \wedge e_k)^{\otimes( m_k-m_{k+1})}\otimes 
(e_1\wedge\dots\wedge e_n)^{\otimes m_n}.
$$
The matrix element 
$$
\la \rho_\bfm(g)\xi_\bfm,\,\xi_\bfm\ra=\prod_{k=1}^{n-1} \Delta_k^{m_k-m_{k+1}}\cdot \det(g)^{m_n}.
$$
 is polynomial if and only if $m_n\ge 0$. 
On the other hand, for $m_n\ge 0$ the representation $\rho_\bfm$ is a polynomial representation
 by the construction.

\section{Proofs}

\COUNTERS

{\bf \punct Direct sums.%
\label{ss:sums-1}} {\sc Proof of Theorem \ref{th:sums}.a.}  Take an element
$g\in \U(\alpha+mi)$ written as
\begin{equation}
g=\left(\begin{array}{c|ccc}
a&b_1&b_2&\dots\\
\hline
c_1&d_{11}&d_{12}&\dots\\
c_2&d_{21}&d_{22}&\dots\\
\vdots&\vdots&\vdots&\ddots
\end{array}
\right)
\label{eq:g}
\end{equation}
and an element $\wt g\in \U(\beta+mj)$ written as
\begin{equation}
\wt g=\left(\begin{array}{c|ccc}
\wt a&\wt b_1&\wt b_2&\dots\\
\hline
\wt c_1&\wt d_{11}&\wt d_{12}\vphantom{\wt {A^{A^A}}}&\dots \\
\wt c_2&\wt d_{21}&\wt d_{22}&\dots\\
\vdots&\vdots&\vdots&\ddots
\end{array}
\right)
\label{eq:g-gtilde}
\end{equation}
Consider the  block matrix
of order
$$\alpha+\beta+i+j+\dots+i+j=(\alpha+\beta)+\underbrace{(i+j)+\dots+(i+j)}_{\text{$m$ times}}
$$
given by
$$
g(\oplus)\wt g:=
\left(
\begin{array}{cc|ccccc}
a&0&b_1&0&b_2&0&\dots\\
0&\wt a&0&\wt b_1&0&\wt b_2&\dots\\
\hline
c_1&0&d_{11}&0&d_{12}&0&\dots\\
0&\wt c_1&0&\wt d_{11}&0&\wt d_{12}&\dots\\
c_2&0&d_{21}&0&d_{22}&0&\dots\\
0&\wt c_2&0&\wt d_{21}&0&\wt d_{22}&\dots\\
\vdots&\vdots&\vdots&\vdots&\vdots&\vdots&\ddots
\end{array}
\right).
$$
By formula \eqref{eq:chi-sigma},
\begin{equation*}
\Theta[g(\oplus)\wt g;S]=\Theta[g;S]\oplus\Theta[\wt g;S]
\end{equation*}
(we applied the Krein--Shmul'yan map determined by the matrix
$g(\oplus)\wt g$ to a matrix $S\otimes 1_{i+j}$).

\sm 

{\bf \punct Pointwise products.%
\label{ss:products}} Theorem \ref{th:products} was observed in \cite{Ner-AAMP}.
To be complete we present here the corresponding operation on colligations.
Let $g\in \U(\alpha+mi)$ be represented as \eqref{eq:g} and $\wt g\in \U(\alpha+mj)$
as \eqref{eq:g-gtilde}. We define the matrix $g\odot\wt g$ by
\begin{multline}
g\odot\wt g:=\\ \left(\begin{array}{c|ccccc}
a&b_1&0&b_2&0&\dots\\
\hline
c_1&d_{11}&0&d_{12}&0&\dots\\
0&0&1_j&0&0&\dots\\
c_2&d_{21}&0&d_{22}&0&\dots\\
0&0&0&0&1_j&\dots\\
\vdots&\vdots&\vdots&\vdots&\vdots&\ddots
\end{array}
\right)
\left(\begin{array}{c|ccccc}
\wt a&0&\wt b_1&0&\wt b_2&\dots\\
\hline
0&1_i&0&0&0&\dots\\
\wt c_1&0&\wt d_{11}&0&\wt d_{12}&\dots \\
0&0&0&1_i&0&\dots\\
\wt c_2&0&\wt d_{21}&0&\wt d_{22}&\dots\\
\vdots&\vdots&\vdots&\vdots&\vdots&\ddots
\end{array}
\right)
\end{multline}
Then
$$
\Theta[g\odot \wt g;S]=\Theta[g;S]\,\Theta[\wt g;S].
$$

{\bf\punct Pointwise tensor products.%
\label{ss:tensors}}
{\sc Proof of Theorem \ref{th:tensors}.} This is a corollary of the previous statement.
Since 
$$
A\otimes B=(1\otimes B)\cdot(A\otimes 1),
$$
it is suffictient to verify the claim for $F_1\otimes 1$ and $1\otimes F_2$.

The function $F_1\otimes 1$ is a direct sum of several copies of the function $F_1$.
By Theorem \ref{th:sums} it is a caracteristic function.
More presisely, if $F_1$ is generated by a matrix 
$\begin{pmatrix}p&q\\r&t\end{pmatrix}$,
then
\begin{multline}
\Theta\left[\begin{pmatrix}p&q_1&\dots&q_m\\
r_1&t_{11}&\dots&t_{1m}\\
\vdots&\vdots&\ddots&\vdots\\
r_m&t_{m1}&\dots&t_{mm}
\end{pmatrix};S \right]\otimes 1_\beta=\\=
\Theta\left[\begin{pmatrix}
p\otimes 1_\beta&q_1\otimes 1_\beta&\dots&q_m\otimes 1_\beta\\
r_1\otimes 1_\beta&t_{11}\otimes 1_\beta&\dots&t_{1m}\otimes 1_\beta\\
\vdots&\vdots&\ddots&\vdots\\
r_m\otimes 1_\beta&t_{m1}\otimes 1_\beta&\dots&t_{mm}\otimes 1_\beta
\end{pmatrix};S \right]
\label{eq:chi-otimes-1}
\end{multline}

The identity
\begin{equation}
1_\alpha\otimes\Theta\left[\begin{pmatrix}a&b\\c&d\end{pmatrix};S \right] =
\Theta\left[\begin{pmatrix}1_\alpha\otimes a&1_\alpha\otimes b\\ 1_\alpha\otimes c&1_\alpha\otimes d\end{pmatrix};S \right]
\label{eq:chi-1-otimes}
\end{equation}
immediately follows from \eqref{eq:chi-sigma}. Therefore
$$
\Theta\left[\begin{pmatrix}
p&q\\r&t
\end{pmatrix};S \right]\otimes \Theta\left[\begin{pmatrix}
a&b\\c&d
\end{pmatrix};S \right]
$$
is a characteristic function of the colligation
\begin{equation}
\left(
\begin{array}{c|ccccc}
p\otimes a&q_1\otimes 1_\beta&p\otimes b_1& q_2\otimes 1_\beta&p\otimes b_2&\dots\\
\hline
r_1\otimes a& t_{11}\otimes 1_\beta& r_1\otimes b_1& t_{12}\otimes 1_\beta& r_1\otimes b_2&\dots\\
1_\alpha\otimes c_1&0&1_\alpha\otimes d_{11}& 0&1_\alpha\otimes d_{12}&\dots\\
r_2\otimes a& t_{21}\otimes 1_\beta& r_2\otimes b_1& t_{22}\otimes 1_\beta& r_2\otimes b_2&\dots\\
1_\alpha\otimes c_2&0&1_\alpha\otimes d_{21}& 0&1_\alpha\otimes d_{22}&\dots\\
\vdots& \vdots&\vdots&\vdots&\vdots&\ddots
\end{array}
\right).
\end{equation}

\sm

{\bf \punct Compositions.%
\label{ss:compositions}} {\sc Proof of Theorem \ref{th:compositions}.}
By \eqref{eq:chi-sigma} and \eqref{eq:chi-1-otimes} we have
\begin{multline*}
G\circ F(S)=\sigma\left[\begin{pmatrix}a&b\\c&d\end{pmatrix};
1_j\otimes \sigma\left[\begin{pmatrix}
p&q\\r&t
\end{pmatrix};1_i\otimes S \right]\right]=\\=
G\circ F(S)=\sigma\left[\begin{pmatrix}a&b\\c&d\end{pmatrix};
\sigma\left[\begin{pmatrix}
1_j\otimes p&1_j\otimes q\\1_j\otimes r&1_j\otimes t
\end{pmatrix};1_i\otimes S \right]\right].
\end{multline*}
If $\det\bigl(1_{\beta j}-d(1_j\otimes p)\bigr)\ne 0$ (condition \eqref{eq:cond1}), then we can apply
formula \eqref{eq:circledast} and Lemma \ref{l:krein-schm}. In this case we come to
$$
\sigma\left[\begin{pmatrix}a&b\\c&d\end{pmatrix} \circledast \begin{pmatrix}
1_j\otimes p&1_j\otimes q\\1_j\otimes r&1_j\otimes t
\end{pmatrix}; 1_{ij}\otimes S  \right],
$$
where $\circledast$-product of matrices is defined by \eqref{eq:circledast}.
Since the both 'factors' are unitary matrices,  we get an inner map $\B_\alpha\to \B_\gamma$.

\begin{corollary}
	\label{c:aut} Let $F\in \Char[m,\alpha]$.
	
	\sm 
	
	{\rm a)} For $h\in \U(\alpha, \alpha)$ we have $\gamma[h]\circ F\in \Char[m,\alpha]$.
	
	\sm 
	
	{\rm b)} For $h'\in \U(m,m)$ we have $F\circ \gamma[h']\in \Char[m,\alpha]$.
\end{corollary}

Indeed, in these cases our condition holds.
\hfill $\square$

\sm 

Let us finish the proof of Theorem \ref{th:compositions}.
Let condition \eqref{eq:cond2} holds, i.e., $\det\bigl(1_{\beta j}-d\,(1_j\otimes F(S_0)) \bigr)\ne 0$. We take $h\in \U(\alpha,\alpha)$ sending 0 to $S_0$, then $F(\gamma[h;0])=F(S_0)$,
  we  apply Corollary \ref{c:aut} to $G\circ(F\circ \gamma[h])$ and refer to the already
  proved part of Theorem \ref{th:compositions}. Thus $G\circ (F\circ \gamma[h])\in\Char(\alpha,\gamma)$.
  Next, we apply Corollary \ref{c:aut} to  $(G\circ F\circ \gamma[h])\circ \gamma[h^{-1}]$ and get the desired statement.
  
  \sm
  
  {\bf \punct Splitting off summands.%
  \label{ss:sums-2}} {\sc Proof of Theorem \ref{th:sums}.b.}
  Let a characteristic function $F\in\Char[m,\alpha+\beta]$
  have a block  form
  $$
  F(z):=\begin{pmatrix}
  F_1(z)&0\\
  0&F_2(z)
  \end{pmatrix}.
  $$
  Let us show that $F_1\in \Char[m,\alpha]$.
  
  \sm 

First, assume that $F_2\in \Inn_\circ[m,\beta]$. 
We consider a Krein--Shmul'yan map $G:\B_{\alpha+\beta}\to \B_\alpha$
 determined by the matrix
 \begin{equation}
 \left(
 \begin{array}{c|cc}
 0&1_\alpha&0\\
 \hline 
 1_\alpha &0&0\\
 0&0&1_\beta
 \end{array}
 \right)
 \label{eq:liv}
 \end{equation}
 and take the composition $G\circ F$,
 \begin{multline*}
 G\circ F(z)=\\=\begin{pmatrix}1_\alpha&0 \end{pmatrix}
 \begin{pmatrix}F_1(z)&0\\0&F_2(z)\end{pmatrix}
 \left\{\begin{pmatrix}1_\alpha&0\\0&1_\beta\end{pmatrix} -
 \begin{pmatrix}0&0\\0&1_\beta \end{pmatrix}
 \begin{pmatrix}F_1(z)&0\\0&F_2(z)\end{pmatrix}
  \right\}^{-1}
  \begin{pmatrix}1_\alpha\\0 \end{pmatrix}.
 \end{multline*}
 The matrix in curly brackets is
 $$
 \begin{pmatrix}
 1_\alpha&0\\
 0&1_\beta-F_2(z)
 \end{pmatrix},
 $$
 it is invertible, and by Theorem \ref{th:compositions} the map $G\circ F$ is contained in
 $\Char[m,\alpha]$. But 
 $$G\circ F=F_1(z)$$
 and this implies our statement.

 \sm
 
 Second, let $F_2\notin \Inn_\circ(m,\alpha)$. Then $F_2(\B_m)$ is contained in some component
 $C$ of the boundary of $\B_\beta$. By Corollary \ref{c:aut} we can assume that $C$ is in a canonical
 form, i.e., $C$ consists of matrices 
 $$
 \begin{pmatrix}
 u&0\\0&1_{k}
 \end{pmatrix}, \qquad \text{where $u\in \B_{\beta-k}$.}
 $$
 and therefore $F_2(z)$ has the form
 \begin{equation}
 F_2(z)=\begin{pmatrix}
 R(z)&0\\0&1_{k}
 \end{pmatrix},\text{where $\|R(z)\|<1$ for $z\in \B_m$.}
 \label{eq:FR}
 \end{equation}
 Now we choose $\lambda\in\C$ such that $|\lambda|=1$ and $\lambda\ne1$.
 Instead of \eqref{eq:liv}
 we take the matrix
 $$
 \left(
 \begin{array}{c|cc}
 0&1_\alpha&0\\
 \hline 
 1_\alpha &0&0\\
 0&0&\lambda\cdot 1_\beta
 \end{array}
 \right)
 $$
and the corresponding Krein--Shmul'yan map $G$. By \eqref{eq:FR} the matrix  
$$
\left\{\begin{pmatrix}1_\alpha&0\\0&1_\beta\end{pmatrix} -
\begin{pmatrix}0&0\\0&\lambda 1_\beta \end{pmatrix}
\begin{pmatrix}F_1(z)&0\\0&F_2(z)\end{pmatrix}
\right\}=\begin{pmatrix}
1_\alpha&0\\
0&1_\beta- \lambda F_2(z)
\end{pmatrix}
$$
is invertible, and by Theorem \ref{th:compositions} we have
$G\circ F=F_1(z)\in \Char(m,\alpha)$.

\sm

{\bf \punct Compositions with polynomial representations.%
\label{ss:tensors-over}} {\sc Proof of Theorem \ref{th:representations}.}
By Theorem \ref{th:sums}.a, it is sufficient to consider irreducible representations.
By construction given in Subsect. \ref{ss:representation}, any irreducible polynomial representation
$\rho_\bfm=\rho_{m_1,\dots,m_n}$ is contained in tensors
$$
\bigotimes_{k=1}^n  \Bigl(\bigwedge\nolimits^k \C^n\Bigr)^{\otimes(m_k-m_{k-1})}\,\otimes\,
 \Bigl(\bigwedge\nolimits^n
\C^n
\Bigr)^{m_n}\subset \bigl(\C^n\bigr)^{\otimes\bigl( \sum_{k=1}^n m_k\bigr)}.
$$
By Theorem \ref{th:tensors} for any characteristic function $F(z)$ a function
$F(z)^{\otimes L}$ is a characteristic function.
By Theorem \ref{th:sums}.b we can split off a direct summand.

\sm 

{\bf \punct Boundary components.%
\label{ss:boundary-restriction}} {\sc Proof of Theorem \ref{th:boundary}.a.}
Without loss of generality we can assume that $C$ has the canonical form
$\begin{pmatrix}
u&0\\0&1_k
\end{pmatrix}$. Therefore our function $F$ splits into a direct sum. By Theorem \ref{th:sums}.b
we can split off a summand.

\sm 

{\sc Proof of Theorem \ref{th:boundary}.b}
Again, we can assume that $C\subset \B_m$ consists of matrices
$\begin{pmatrix}
u&0\\0&1_l
\end{pmatrix}
$. As in Subsect. \ref{ss:compositions} we can assume that $S_0:=\begin{pmatrix}
0&0\\0&1_l
\end{pmatrix}$. The identical embedding $u\mapsto \begin{pmatrix}
u&0\\0&1_l
\end{pmatrix}$ is  a Krein--Shnul'yan map determined by the matrix
$$
\left(
\begin{array}{cc|c}
0&0& 1_{m-l}\\
0&1_l&0\\
\hline
1_{m-l}&0&0
\end{array}
\right).
$$
Now we can apply Theorem \ref{th:compositions}.

\tt
\noindent
Math. Dept., University of Vienna\\
\&Institute for Theoretical and Experimental Physics (Moscow); \\
\&MechMath Dept., Moscow State University;\\
\&Institute for Information Transmission Problems;\\
URL: http://mat.univie.ac.at/$\sim$neretin/

\end{document}